\documentclass[,english]{smfart}
\usepackage{MnSymbol}
\usepackage{smfthm}
\usepackage[matrix, arrow,all,cmtip,color]{xy}
\usepackage[utf8]{inputenc}
\usepackage{epigraph}
\usepackage{url} 
\usepackage{hyperref}

\def\VV{\mathcal Var}
\def\Var{\mathcal Var}
\def\VVK{\VV_K}
\def\FF{{\mathcal F}}
\def\Points{{\rm Ob}\,}
\def\Aut{Aut}
\def\Paths{{\rm Mor}\,}
\def\Mor{{\rm Mor}\,}
\def\Mot{{\mathcal Mot}}

\def\GSpZ{{\rm GSp_\Z}}
\def\C{{\Bbb C}}
\def\Q{{\Bbb Q}}
\def\Qbar{{\bar\Q}}
\def\kbar{{\bar{k}}}

\def\Z{{\Bbb Z}}

\def\AA{{\Bbb A}}
\def\Ql{{\Bbb Q}_l}
\def\Zl{{\Bbb Z}_l}
\def\Zhat{{\Bbb {\hat Z}}}

\def\Fpbar{{\Bbb {\bar F}}_p}
\def\Groupoids{Groupoids}

\def\Gal{{\rm Gal}}

\def\xra{\xrightarrow}

\def\lra{\longrightarrow}
\def\Lra{\Longrightarrow}

\def\char{{\rm char}\,}
\def\card{{\rm card}\,}
\def\Ob{{\rm Ob}\,}
\def\gen#1{\left< #1\right>}
\def\genn#1#2{\left< #1 \right>_{#2}}

\def\QVect{\Q{\textrm{-Vect}}}
\def\Qvect{\QVect}

\def\ZVect{\Z{\textrm{-Mod}}}

\def\ZlVect{\Zl{\textrm{-Mod}}}
\def\Zhatvect{\Zhat{\textrm{-Mod}}}

\def\Het{H_{\textrm et}}
\def\Hetl{H_l}
\def\Hl{H_l}
\def\tensor{\otimes}

\def\Het{H_{et}}
\def\Ql{{\Bbb Q}_l}

\def\Lwlw{L_{\omega_1\omega}}
\def\Hsing{H_{\text{sing}}}
\def\Q{{\Bbb Q}}
\def\Qbar{\bar{\Bbb Q}}

\def\lra{\longrightarrow}
\def\Lra{\Longrightarrow}
\def\xra{\xrightarrow}
\def\xla{\xleftarrow}

\def\R{\Bbb R}

\def\Qvect{\Q\text{-Vect}}
\def\QVect{\Q\text{-Vect}}
\def\QHodge{\Q\text{-Hodge}}

\def\QpVect{\Q_p\text{-Vect}}

\def\id{id}

\def\Var{\mathcal Var}
\def\Mot{\mathcal Mot}
\def\V{\mathcal V}
\def\VV{\mathcal V}
\def\WW{\mathcal WW}

\theoremstyle{plain}

\newtheorem{problem}{Problem}
\newtheorem{defn}{Definition}
\newtheorem{example}{Example}
\newtheorem{question}{Question}

\theoremstyle{definition}
\newtheorem{remark}{Remark}

\def\qftp{{\mathrm{qftp}}} 
\newcommand{\xdashleftrightarrow}[2][]{\ext@arrow 3359\leftrightarrowfill@@{#1}{#2}}

\def\ex{\mathrm{ex}}
\def\Ker{\mathrm{Ker}\,}

\usepackage{color}

\newcommand\longdashto{\mathrel{
  -\mkern-6mu{\longrightarrow}\mkern-28mu{\color{white}\cdots}\mkern12mu
}}

\title[Standard conjectures in model theory]{
standard conjectures in model theory, and categoricity of comparison isomorphisms
}
\author[misha gavrilovich]{notes by misha gavrilovich
}

\address{National Research University Higher School of Economics, Saint-Petersburg\\
Institute for Problems of Regional Economics RAS                                                                       
38 Serpuhovskaya st., Saint-Petersburg }
\email{mi\!\!\!ishap\!\!\!p@sd\!\!\!df.org}
\urladdr{http://mishap.sdf.org/hcats.pdf}
\thanks{Some ideas were formed and formulated with help of Martin Bays.  
$\tt{mi\!\!\!ishap\!\!\!p@sd\!\!\!df.org}$. \href{http://mishap.sdf.org/hcats.pdf}{\tt http://mishap.sdf.org/hcats.pdf}
}

\begin{document}

\maketitle

\begin{abstract}
We formulate two conjectures about \'etale cohomology and fundamental groups motivated by categoricity conjectures in model theory. 

One conjecture says that there is a unique $\mathbb Z$-form
of the \'etale cohomology of complex algebraic varieties, up to $Aut(\mathbb C)$-action
on the source category;
put differently, each comparison isomorphism
between Betti and \'etale cohomology comes from a choice of a topology on $\mathbb C$.

Another conjecture says that each functor to groupoids from the category of complex algebraic varieties
which is similar to the topological fundamental groupoid functor $\pi_1^{top}$, 
in fact factors through $\pi_1^{top}$, up to a field automorphism of the complex numbers
acting on the category of complex algebraic varieties. 

We also try to present some evidence towards these conjectures, 
and show that some special cases seem related to Grothendieck 
standard conjectures and conjectures about motivic Galois group.  
\end{abstract}

\section{Introduction}

We consider the following question as it would be understood by a model theorist 
\begin{enonce*}{Question} Is there a purely algebraic definition of the notion of 
singular (Betti) cohomology or the topological fundamental groupoid of a complex algebraic variety?
\end{enonce*}
and formulate precise conjectures proposing that comparison isomorphism of \'etale cohomology/fundamental groupoid 
admits such a purely algebraic definition (characterisation). These conjectures are direct analogues of categoricity 
theorems and conjectures in model theory, particularly those on pseudoexponentiaton [Zilber].

We then show that 
some special cases of these conjectures seem related to Grothendieck
standard conjectures and conjectures about motivic Galois group, particularly the image of $l$-adic Galois representations.

Note that an algebraic geometer might interpret 
 the question differently and in that interpretation, the answer is well-known to be negative. 

We now explain our motivation in two essentially independent ways. \S1.1 explains how a model theorist 
would interpret the question above; \S1.2 views these conjectures as continuation of work in model theory on 
the complex field with pseudoexponentiation [Zilber, Bays-Zilber, Bays-Kirby. Manin-Zilber]
and its main goal is to make the reader aware of the possibilities offered by methods of model theory.

\subsection{How to interpret the question.}
Let us now explain the difference between how an algebraic geometer and a model theorist might interpret the question.

Let $H_{\mathrm{top}}$ be a functor defined on the category $\Var$ of algebraic varieties (say, separated schemes of finite type) 
over the field $\C$ of complex numbers; we identify this category with a subcategory of the category of topological spaces. 
We shall be interested in the case when $H_{\mathrm{top}}$ is either the functor $\Hsing:\Var\lra Ab$ of singular cohomology
or the fundamental groupoid functor $\pi_1^{top}:\Var\lra Groupoids$. 

An algebraic geometer might reason as follows. 
A purely algebraic definition applies both to $H_{\mathrm{top}}$ and $H_{\mathrm{top}}\circ \sigma$
where $\sigma:\C\lra\C$ is a field automorphism. Hence, to answer the question in the negative,
it is enough to find a field automorphism $\sigma:\C\lra\C$ such that $H_{\mathrm{top}}$ and $H_{\mathrm{top}}\circ \sigma$ 
differ. And indeed, [Serre, Exemple] constructs an example of a projective algebraic variety $X$ and a field automorphism $\sigma$
such that $X(\C)$ and $X^\sigma(\C)$ have non-isomorphic fundamental groups. 

A model theorist might reason as follows. A purely algebraic definition applies both to $H_{\mathrm{top}}$ and $H_{\mathrm{top}}\circ \sigma$
where $\sigma:\C\lra\C$ is a field automorphism. Hence, we should try to find purely algebraic 
description  (possibly involving extra structure) of $H_{\mathrm{top}}$  
which fits precisely functors of form $H_{\mathrm{top}}\circ \sigma, \sigma\in \Aut(\C)$ with the extra structure. 
We say that {\em such a purely algebraic description describes $H_{\mathrm{top}}$ (with the extra structure) 
uniquely up to an automorphism of $\C$.}

For $H_{\mathrm{top}}=\Hsing$ the singular (Betti) cohomology theory, a model theorist might continue thinking as follows.
The singular (Betti) cohomology theory admits a comparison isomorphism to a cohomology theory defined purely
algebraically, say $l$-adic \'etale cohomology theory. This is an algebraic description in itself. 
However, note that it considers the $l$-adic \'etale cohomology theory and the comparison isomorphism 
as part of structure. 
Thus an appropriate conjecture (see~\S\ref{htop:uniq}) 
gives a purely algebraic description of the family of comparison isomorphisms
coming from a choice of topology on $\C$
$$
\Hsing(X(\C)_{top},\Z)\tensor \Z_l \xra{\ \sigma\ } \Het(X(\C),\Zl),\ \ \ \sigma\in\Aut(\C)$$

For $H_{\mathrm{top}}=\pi_1^{top}$ a model theorist might continue thinking as follows. 
The profinite completion of the topological fundamental groupoid functor is the \'etale fundamental groupoid 
defined algebraically. This is an algebraic property of the topological fundamental groupoid
on which we can base our purely algebraic description 
if we include the \'etale fundamental groupoid as part of structure.
Essentially, this describes subgroupoids of the \'etale fundamental groupoids. 
Category theory suggests to consider a related universality property (see~\S\ref{pitop:uniq}):
{\em up to $\Aut(\C)$ action on the category of complex algebraic varieties},
there is a universal functor among those
 whose profinite
completion embeds into the \'etale fundamental groupoid,
and it is the topological fundamental groupoid.
Some technicalities may be necessary to ignore non-residually finite fundamental groups.
\subsection{Pseudo-exponentiation, Schanuel conjecture and categoricity theorems in model theory.}
Complex topology allows to construct a number of objects with good algebraic properties 
e.g.~a group homomorphism 
$\exp:\C^+\lra \C^*$, singular (Betti) cohomology theory  
and the topological fundamental groupoid of varieties of complex algebraic varieties. 

A number of theorems and conjectures says that such an object constructed topologically or analytically is ``free'' or ``generic'', for lack of better term, in the sense that it  satisfies algebraic relations only, or mostly, for ``obvious'' reasons 
of algebraic nature. 

Sometimes such a conjecture is made precise by saying that a certain
automorphism group is as large as possible subject to some ``obvious
obstructions or relations'' imposed by functoriality and/or homotopy theory.
Such an automorphism group may involve values of functions
 or spaces defined analytically or topologically. 

A natural question to ask is whether these conjectures are ``consistent" in
the sense that there do exist such ``free" objects with the conjectured
properties, not necessarily of analytic or topological origin. 

Methods of model theory allow to build such objects by an elaborate transfinite induction. 
In what follows we shall sketch results of [Zilber, Bays-Kirby] which does this for
the complex exponential function and Schanuel conjecture. 

Let us now explain what we mean by showing how to view Kummer theory, Hodge
conjecture, conjectural theory of the motivic Galois group, and Schanuel
conjecture in this way.

\subsubsection{Kummer theory.}
An ``obvious way'' to make 
 $e^{\alpha_1/N},...,e^{\alpha_n/N}$, $N>0$
satisfy a polynomial relation is 
is to pick ${\alpha_1},...,{\alpha_n}$ 
such that they satisfy a $\Q$-linear relation over $2\pi i$,
which is preserved by $\exp$, or 
such that that $e^{\alpha_1/M},...,e^{\alpha_n/M}$
satisfy a polynomial relation for some other $M$. 

Kummer theory tells you these are the only reasons 
for polynomial relations between these numbers. 
This is stated precisely in terms of automorphisms groups as follows:

For any $\Q$-linearly independent numbers ${\alpha_1},...,{\alpha_n}\in \C$
there is $N>0$ such that for any $m>0$ it holds
$$\Gal(\Q(e^{\frac{\alpha_1}{mN}},...,e^{\frac{\alpha_n}{mN}},e^{2\pi i\Q})/
\Q(e^{\frac{\alpha_1}N},...,e^{\frac{\alpha_n}N},e^{2\pi i\Q}))\approx (\Z/m\Z)^n$$

\subsubsection{Hodge conjecture}
Consider the Hodge theory of a non-singular complex projective manifold $X(\C)$. 
By Chow theory we know that $X$ is in fact a complex algebraic variety
and an easy argument using harmonic forms shows that an algebraic subvariety $Z(\C)$ defines
an element of 
$H(X,\Q)\cap H^{(p,p)}(X,\C)$
where $H^{(p,p)}(X,\C)$ is a certain linear subspace of $H^{2p}(X,\C)$ defined analytically.

A topological cycle in $X(\C)$ defines an element of $H(X,\Q)$ which may lie in $H^{(p,p)}(X,\C)$.
An ``obvious reason'' for this is that it comes from an algebraic subvariety, 
or a $\Q$-linear combination of such. Hodge conjecture tells you that this is the only reason it could happen.

\subsubsection{$l$-adic Galois representations and motivic Galois group $\Aut^\tensor({\Hsing}_\sigma)$}
Remarks below are quite vague but we hope some readers might find them helpful.  In \S\ref{serre} 
we sketch several definitions and conjectures in the conjectural theory of motivic Galois group following [Serre].

We would like to think that these conjectures say that the singular (Betti) cohomology theory of complex algebraic varieties
is ``free'' in the sense
that it satisfies algebraic relations only, or mostly, for ``obvious'' reasons of algebraic
nature. The theory of the motivic Galois group assumes that there are many automorphisms of 
the singular cohomology theory of complex algebraic varieties, and they form a
pro-algebraic, in fact pro-reductive ([Serre, Conjecture~2.1?], group.  
Conjectures on $l$-adic Galois representations, e.g.~[Serre, Conjecture~3.2?,9.1?] 
describe the image of Galois action as being dense or open in a certain 
algebraic group defined by cohomology classes which Galois action has to preserve (or is conjectured to preserve). 

Let us very briefly sketch some details. 

The conjectural theory of the motivic Galois group [Serre], also cf.~\S\ref{serre},
assumes that the following is a well-defined algebraic group:  
$$G_E=\Aut^\tensor({\Hsing}_\sigma:\left<E\right>\lra\QVect)$$
Here $\sigma:k\lra \C$ is an embedding of a number field $k$ into the field of complex numbers,
$E$ is a pure motive in the conjectural category $\Mot_k$
of pure motives defined over $k$, and 
$\left<E\right>$ is the least Tannakian subcategory of $\Mot_k$ 
containing $E$, and ${\Hsing}_\sigma$ is the fibre functor on $\left<E\right>$
corresponding to the singular cohomology of complex algebraic varieties
and embedding $\sigma:k\lra \C$. This is well-defined if we assume certain conjectures, 
e.g.~Standard Conjectures and Hodge conjecture [Serre, Grothendieck, Kleiman].

[Serre, Conjecture 3.1?] says that $G_E$ is the subgroup of
$GL({\Hsing}_\sigma(E))$  
preserving the tensors corresponding to morphisms
 ${\boldsymbol 1}\lra E^{\tensor r} \tensor E^{\vee\tensor s}$, $r,s\geq 0$. 
Think of these tensors as ``obvious relations'' which have to be preserved.

[Serre, Conjecture 3.2? and Conjecture 9.1?] describe the 
image of $l$-adic Galois representations in $G_E(\Ql)$.

Both say it is dense or open in the group of $l$-adic points of 
a certain algebraic subgroup of $GL_N$; we think of this subgroup 
as capturing ``obvious obstructions or relations'' imposed by functoriality 
of $\Hsing$.

\subsubsection{Schanuel conjecture: questions}

Schanuel conjecture says that for 
$\Q$-linearly independent $x_1,...,x_n\in \C$,
the transcendence degree of $x_1,...,x_n,e^{x_1},...,e^{x_n}$ is at least $n$:
$$\mathrm{tr.deg.}_\Q(x_1,...,x_n,e^{x_1},...,e^{x_n})\geq \mathrm{lin.deg.}_{\Q}(x_1,...,x_n) \eqno{\mathrm{(SC)}}$$
The bound becomes sharp if we use surjectivity to pick $x_2=e^{x_1}$, ..., $x_{i+1}=e^{x_i},...,x_n=e^{x_{n-1}}$
and $e^{x_n}\in \Q$:
$$\mathrm{tr.deg.}_\Q(x_1,e^{x_1},e^{e^{x_1}},...,x_n,e^{x_1},e^{e^{x_1}},...,e^{x_n})\leq n$$

Here ``an algebraic relation'' is a polynomial relation between $x_1,...,x_n,e^{x_1},...,e^{x_n}$;
an obvious way to make these numbers satisfy such a relation 
is  to pick $x_i$ such that either $x_i=a_1x_1+...+a_{i-1}x_{i-1}$ 
or $e^{x_i}=a_1x_1+...+a_{i-1}x_{i-1}$ or $x_i=e^{a_1x_1+...+a_{i-1}x_{i-1}}$ 
 where $a_1,...,a_i\in \Q$ are rational.


Is Schanuel conjecture ``consistent'' in the sense that there is a {\em pseudo-exponentiation}, 
i.e.~a group homomorphism $\ex:\C^+\lra \C^*$ satisfying conjectural properties of 
complex exponentiation, in particular Schanuel conjecture?
Does there exist such a ``free''  pseudo-exponentiation  $\ex:\C^+\lra \C^*$, 
e.g.~such that a system of  exponential-polynomial equations has a zero only iff it 
does not contradict Schanuel conjecture? Can we build such an algebraic ``free'' object
without recourse to topology? 

Does every such ``free'' object come from a
choice of topology on $\C$, i.e.~is 
the complex exponential $\exp:\C^+\lra\C^*$ up an automorphism of $\C$? 

Note that the last question  is the only one which mentions topology. 
It turns out this difference is crucial: model theory says nothing about this question
while giving fairly satisfactory positive answers to the previous ones.

\subsubsection{Schanuel conjecture and pseudoexponentiation: answers} 
The following theorem of [Zilber]  provides a positive answer for $\exp:\C^+\lra\C^*$. 
 For a discussion of the theorem and surrounding model theory see [Manin-Zilber, 6.16];
for a proof, detailed statements and generalisations to other analytic functions
see [Bays-Kirby, Thm.~1.2,Thm.~1.6; Thm.~9.1; also Thm.~8.2; Thm.~9.3] and references therein.
\def\ex{\mathrm{ex}}
\def\Ker{\mathrm{Ker}\,}
\begin{enonce}{Theorem}[Zilber]\label{Zil:preudoexp}
 Let $K$ be an uncountable algebraically closed field of characteristic $0$.

Up to $\Aut(K)$, there is a unique surjective group homomorphism 
$\ex:K^+\lra K^*$ 
\begin{itemize}
\item[(SK)] (Standard Kernel) $\Ker \ex$ is the infinite cyclic group generated by a transcendental element
\item[(SC)] (Schanuel Property) $\ex:K^+\lra K^*$ satisfies Schanuel conjecture
\item[(SEAC)](Strong exponential-algebraic closedness)   
 any system of $n$ independent exponential-polynomial equations in $n$
     variables that does not directly contradict Schanuel conjecture has a regular
     zero, but not more than countably many
\end{itemize}
\end{enonce}

Call this unique group homomorphism $\ex:K^+\lra K^*$ {\em pseudoexponentiation} defined on field $K$.

In somewhat more detail, this can also be expressed as follows.

Let $K$ and $K'$ be two uncountable algebraically closed fields of characteristic $0$, 
and let $\ex:K^+\lra K^*$ and $\ex':K'^+\lra K'^*$ be group homomorphisms satisfying the properties above.

Then if there is a bijection $\sigma_0:K\xra{\,\,\approx\,\,} K'$, then there is a bijection $\sigma:K\xra{\,\,\approx\,\,} K'$
preserving $+$, $\cdot$, and $\ex$, i.e.~such that for each $x,y\in K$ it holds
$$\sigma(x+y)=\sigma(x)+\sigma(y),\ \ \sigma(xy)=\sigma(x)\sigma(y),\ \ 
\sigma(\ex(x))=\ex'(\sigma(x))$$

\begin{enonce}{Conjecture}[Zilber] If $\card K=\card \C$, 
then $(K,+,\cdot,\ex)$ is isomorphic to $(\C,+,\cdot,\exp)$.
\end{enonce}

Our conjectures are direct analogues of the Theorem and Conjecture above
stated in the language of functors. Instead of the complex exponentiation
we consider the comparison isomorphisms between topological and \'etale 
cohomology, resp. fundamental groupoid functor. 
We hope that model theoretic  methods used by [Zilber] may be of use 
in proving these conjectures.

\subsubsection{Pseudoexponentiation: automorphisms groups} It is known that certain automorphisms groups 
associated with pseudoexp are largest possible in the following sense. 

We need some preliminary definitions. We say that tuples $a$ and $b$ in $K$ have {\em
the same quantifier-free type}, write $\qftp(a)=\qftp(b)$, 
iff they satisfy the same exponential-polynomial equations, and, moreover, 
the same exponential-polynomial equations with coefficients with $a$, resp. $b$, 
have a solution; see  [Bays-Kirby, \S6, Def.~6.7] for details.
Note that for a finite tuple $a$ in $K$, 
there is a minimal $\Q$-linear vector subpace $A\supset a$ such that $A\leq_\delta K$
and this $A$ determines $\qftp(a)$ (see below for the definition of $\leq_\delta$).

We quote from [Bays-Kirby, Def.~6.1, Proposition~6.5].

\def\qftp{{\mathrm{qftp}}}
\begin{enonce}{Fact} Let $K$ be a field with pseudoexponentiation as defined above.  
\begin{itemize}
\item[QM4.] (Uniqueness of the generic type) Suppose that $C, C ' \subset M$ are countable
     closed subsets, enumerated such that $\qftp(C) = \qftp(C ' )$. If $a \in M\setminus C$
     and $a' \in M\setminus C '$ then $\qftp(C, a) = \qftp(C ' , a' )$ (with respect to the same
     enumerations for $C$ and $C '$ ).
\item[QM5.]  ($\aleph_0$-homogeneity over closed sets and the empty set)
Let $C, C' \subset K$ be countable closed subsets or empty, enumerated such
that $\qftp(C) = \qftp(C' )$, and let $b, b'$ be finite tuples from $K$ such that
$qftp(C, b) = \qftp(C ' , b' )$, and let $a \in cl(C, b)$. Then there is $a' \in  K$ such
that $qftp(C, b, a) = \qftp(C ' , b' , a' )$.

\item[QM5a.] ($\aleph_0$-homogeneity over the empty set)
      If $a$ and $b$ are finite tuples from $K$ and $\qftp(a) = \qftp(b)$ then there is 
a field automorphism $\theta:K\lra K$ preserving $\ex:K^+\lra K^*$ 
such that $\theta(a) = b$.
\end{itemize}
\end{enonce}

Note that it is an open problem to construct a non-trivial automorphism of $(\C,+,\cdot,\exp)$.
\def\base{\mathrm{base}}
\subsubsection{Remarks about the proof.} We adapt [Manin-Zilber, 6.11-6.16];
see also [Bays-Kirby] for a detailed exposition in a more general case using different terminology.
Pseudoexponentiation is constructed by an elaborate transfinite induction. 
We start with an algebraically closed field $K_\base\subset K$ and a partial 
 group homomorphism $\ex:K_\base^+\longdashto K_\base^*$ 
and try to extend 
the field and 
the group homomorphism such that 
it is related to the field in as free a
way as possible. 

Informally the freeness condition is described as follows:
\begin{itemize}\item[(Hr)]
the number of independent explicit basic dependencies {\em added} to 
a subset $X\cup \ex(X)$ of $K$ by the new structure is at most 
the dimension of $X\cup \ex(X)$ in the old structure.
\end{itemize}

This is made precise in the following way. 

{\em The new structure} is the group homomorphism $\ex:K^+\lra K^*$; 
{\em explicit basic dependencies in $X\cup \ex(X)$ added by the new structures} 
are defined as 
as equations $\ex(x)=y$ where $x\in X$. For example, for $X=\{x\}$  where 
$\ex(\ex(x))=x$, we do not regard $\ex(\ex(x))=x$ as a explicit basic dependency in $X\cup\ex(X)=\{x,\ex(x)\}$.  

{\em The number of independent basic explicit dependencies} is the $\Q$-linear dimension 
$\mathrm{lin.dim.}_{\Q}(x_1,...,x_n)$; {\em the dimension of  $X$ in the old structure} is its transcendence degree
which is equal to $\mathrm{tr.deg.}(x_1,...,x_n,\ex(x_1),...,\ex(x_n))$.  

With this interpretation, (Hr) becomes Schanuel conjecture (SC). 

Define {\em Hrushovski predimension} $\delta(X):=\mathrm{tr.deg.}(X\cup\ex(X))-\mathrm{lin.dim.}_{\Q}(X)$. 
Say a partial group homomorphism $\ex:K^+\longdashto K^*$ satisfies {\em Hrushovski inequality
with respect to Hrushovski predimension $\delta$} 
iff 
for any finite $X\subset K$ it holds $\delta(X)\geq 0$. 
An extension $(K,\ex_K)\subset (L,\ex_L)$ of fields equipped with partial group homomorphisms 
is {\em strong}, write $K\leq_\delta L$, iff
all dependencies between elements of $K$ occurring in $L$ can be detected
already in $K$, i.e.~for every finite $X \subset K$,
$$  \min\{\,\delta(Y) : Y \text{ finite, } X \subset Y \subset K\,\} 
= \min\{\,\delta(Y ) : Y \text{ finite, } X \subset Y \subset L\,\}$$

We then build a countable algebraically closed field $(K_{\aleph_0},\ex_{K_{\aleph_0}})$ by
taking larger and larger  strong extensions $K_\base\leq_\delta K_1 \leq K_2 \leq_\delta ... $
of finite degree. 
If we do this with enough care, we obtain a countable algebraically closed field $K_{\aleph_0}=\cup K_n$  
and a group homomorphism $\ex_{K_{\aleph_0}}: K_{\aleph_0}^+\lra K_{\aleph_0}^*$ defined everywhere
which satisfies (SC) and other conditions of Theorem~\ref{Zil:preudoexp}.
For details see [Bays-Kirby, \S5] where it is described in terms of taking Fraisse limit
along a category of strong extensions.

Building an uncountable model requires deep model theory; see 
[Bays-Kirby, \S6] and $[\mathrm{BH^2K^2 14}]$. Let us say a couple of words about this. 
In the inductive construction above, being countable is essential: if we start with an uncountable field, 
we can no longer hope to obtain an algebraically closed field after taking union of countably many extensions
of finite degree. 
Very roughly, it turns out that we can construct composites of countable linearly disjoint
algebraically closed fields this way, and this helps to build an uncountable field with pseudoexponentiation 
and prove it is unique in its cardinality.

\subsubsection{Generalisations and Speculations.} 
[Bays-Kirby] generalises the considerations above 
in a number of ways. In particular, they construct pseudo-exponential maps of simple abelian varieties, including
 pseudo-$\wp$-functions for elliptic curve. 
[Proposition 10.1, \S10, ibid.] relates the Schanuel property of these to the Andr\'e-Grothendieck
conjecture on the periods of 1-motives. They 
 suspect that for abelian varieties the predimension inequality $\delta(X)>0$
also follows from the Andr\'e-Grothendieck periods conjecture, but there are
more complications because the Mumford-Tate group plays a role and so have
not been able to verify it. [\S9.2, ibid.] says it is possible to construct a pseudoexponentiation 
incorporating a counterexample to Schanuel conjecture, by suitably modifying the Hrushovski predimention
and thus the inductive assumption (Hr). [\S9.7, also Thm.~1.7, ibid.] considers differential equations. 

We intentionally leave the following speculation vague. 
\begin{enonce}{Speculation} Can one build a pseudo-singular, or pseudo-de Rham cohomology theory, or a pseudo-topological fundamental group functor 
 of complex algebraic varieties, or an algebra of pseudo-periods
which satisfies a number of conjectures 
such as the Standard Corjectures, the conjectural theory of the motivic Galois group,
the conjectures on the image of $l$-adic Galois representaitons,  Andr\'e-Grothendieck periods
conjecture, Mumford-Tate conjecture, etc.? 
\end{enonce}

\subsection{A glossary of terminology in model theory.} 

We give a very quick overview of basic terminology used in model theory. 
See [Tent-Ziegler; Manin-Zilber] for an introduction into model theory.

In logic, a property is called {\em categorical} iff 
any two structures (models) satisfying the property are necessarily isomorphic. 
A {\em structure} or {\em a model} is usually understood as a set $X$  equipped
with names for certain distinguished subsets of its finite Cartesian powers $X^n$, $n>0$, called {\em predicates},
and also equipped with names for certain distinguished functions between its finite Cartesian powers. 
Names of predicates and functions form a {\em language}.
{\em First order formulas in language $L$} is a particular class of formulas which provide
 names for subsets obtained from the
$L$-distinguished subsets by taking finitely many times intersection, union,
completion, and projection onto some of the coordinates; a formula $\varphi(x_1,..,x_n)$ 
defines the subset $\varphi(M^n)$ of $M^n$ 
consisting of tuples satisfying the formula. A {\em theory in language $L$} is a collection of formulas in language $L$.
{\em A model} of a theory $T$ in language $L$ is a structure in language $L$ such that 
for each $\varphi\in T$ $\varphi(M^n)=M^n$ where $n$ is the arity of $\varphi$.

{\em The first order theory of a structure} consists of all possible names (formulas)
 for the subsets $M^n,n\geq 0$, i.e.~formulas $\varphi$ such that $\varphi(M^n)=M^n$.

A {\em categoricity} theorem in model theory usually says that any two 
models of a first order theory of the same uncountable cardinality are necessarily isomorphic,
i.e.~if there is a bijection between (usually assumed uncountable) models $M_1$
and $M_2$ of the theory, then there is a bijection which preserves the distinguished subsets
and functions. A theory is {\em uncountably categorical} iff it has a unique model, up to isomorphism, 
of each uncountable cardinality. 

The {\em type $\mathrm{tp}(a_1,\!...,a_n)=\{\varphi(x_1,\!..,x_n):\varphi(a_1,\!...,a_n)\text{ holds in } M\}$ 
of a tuple $(a_1,..,a_n)\in M^n$} is the collection of all formulas satisfied by 
the tuple $(a_1,..,a_n)$. A {\em type in a theory} is the type of a tuple in a model of the theory. 
The {\em type $\mathrm{tp}(a_1,...,a_n)=\{\,\varphi(x_1,..,x_n):\varphi(a_1,...,a_n)\text{ holds in } M\}$ 
of a tuple $(a_1,..,a_n)\in M^n$ with parameters in subset $A\subset M$} is the collection of 
all formulas with parameters in $A$ satisfied by 
the tuple $(a_1,..,a_n)$. A {\em type in a theory} is the type of a tuple in a model of the theory. 
Informally, the type of a tuple is a syntactic notion playing the role of an orbit of $\Aut^L(M)$ on $M^n$, 
e.g.~in a situation when we do not yet know whether non-trivial automorphisms of $M$ exist.

In an uncountably categorical first order theory with finitely many predicates and functions
the number of types is at most countable, and the number of types with parameters in a subset $A$ 
has cardinality at most $\card A + \aleph_0=\max(\card A,\aleph_0)$.

%
%

\section{\label{htop:uniq}Uniqueness property of comparison isomorphism of singular and \'etale cohomology of a complex algebraic variety}


A {\em $\Z$-form} of a functor $H_l:\VV\lra\ZlVect$ is 
a pair $(H,\tau)$ consisting of a functor  $H:\VV\lra\ZVect$ and an isomorphism 
$$
H\tensor_\Q \Zl \xra{\tau} H_l 
$$
of functors. 

An example of a $\Z$-form we are interested in is given by the comparison isomorphism between 
\'etale cohomology and Betti cohomology, see [SGA 4, XVI, 4.1], also
[Katz,p.23] for the definitions and exact statements. 

Let 
$\Het:Schemes\lra \ZVect$ be the functor of $l$-adic \'etale cohomology, and let $H_{sing}:Top\lra \ZVect$ be 
the functor of singular cohomology. For $X$ a separated $\C$-scheme of finite type there is a canonical {\em comparison isomorphism}  
$$H_{sing}(X(\C),\Z))\tensor \Zl \approx \Het(X,\Zl)
.$$
This defines a $\Z$-form of the functor of $l$-adic \'etale cohomology $\Het(-,\Zl)$ restricted to the category of separated $\C$-schemes of finite type. 

Let $K$ be an algebraically closed field, let $\VVK$ be a category of varieties over $K$. 
A field automorphism $\sigma:K\lra K$ 
acts  $X\mapsto X^\sigma$ on the category $\VVK$ by automorphisms. Moreover, for each variety $X$ defined over $K$,
 a field automorphism $\sigma$ defines an isomorphism $\sigma_X:X\lra X^\sigma$ of schemes (over $\Z$ or $\Z/p\Z$), and hence 
$$\Het(X,\Zl)\xra{\sigma_*}\Het(X^\sigma,\Zl)
.$$

This defines an action of $\Aut(K)$ on the 
$\Z$-forms of $\Hl$:
$$(H,\tau)\longmapsto (H\circ\sigma, \tau\circ\sigma^{-1}_*)$$
$$ 
H(X^\sigma)\xra{\tau\circ\sigma^{-1}_*} \Het(X,\Zl)
.$$

%
%
We conjecture that, up to action of $\Aut(\C)$ defined above,  
the comparison isomorphism between singular and $l$-adic cohomology of 
is the only $\Z$-form  
of the $l$-adic cohomology theory $\Het(-,\Zl)$:

\begin{enonce}{Conjecture}[$Z(H_{sing},\Hl)$] Up to $Aut(\C)$ action, there is a unique  $\Z$-form of the $l$-adic
cohomology theory functor $\Het(-,\Zl)$ restricted to the category 
of separated $\C$-schemes of finite type 
which respects the cycle map and Kunneth decomposition.

In other words, every comparison isomorphism of a $\Z$- and the $l$-adic cohomology theory of 
 separated $\C$-schemes of finite type
is, up to a field automorphism of $\C$,  the standard comparison isomorphism 
$$H_{sing}(X(\C),\Z)\tensor \Zl = \Het(X,\Zl)
$$
\end{enonce}

The conjecture is intended to be too optimistic; it is probably more reasonable to conjecture uniqueness of $\Z$-form of 
the {\em torsion-free} part of the $l$-adic cohomology. 


Assume  Grothendieck Standard Conjectures and that, in particular 
the $l$-adic cohomology theory factors via the category $\Mot_k$ of pure motives
over a field $k$. 
Then, a Weil cohomology theory (cf.~[Kleiman]) corresponds to a tensor fibre functor 
from the category of pure motives, and we may ask how many 
 $\Z$-forms does have the fibre functor corresponding to the $l$-adic cohomology theory.
Moreover, we may formulate a ``local'' version of the conjecture restricting
the functor to a subcategory generated by a single motive.

\begin{enonce}{Conjecture}[$Z(\Het,\genn E k)$] Let $k$ be a number field. 
Assume Grothendieck Standard Conjectures and that, in particular, 
$l$-adic cohomology factors via the category of pure numerical motives $\Mot_k$
over $k$. 

Let $E$ be a motive and let $\gen E$ be the subcategory of $\Mot_k$ 
generated by $E$,~i.e.~the 
least Tannakian subcategory of $\Mot_k$ containing $E$. 
Up to $\Aut(\bar{k}/k)$-action, 
the functor $\Het(-\tensor\Qbar,\Zl):\genn E k \lra\ZlVect$ has at most finitely many  $\Z$-forms.

Moreover, if $E$ has finitely many $\Z$-forms [Serre,10.2?], 
then the functor 
$\Het(-\tensor \Qbar, \Zhat ):\genn E k \lra \Zhatvect$
has at most finitely many  $\Z$-forms.
\end{enonce}

\subsection{An example: an Abelian variety.} Let us give an example of a particular case of the conjecture 
which is easy to prove.

\begin{example} 
Let $A$ be an Abelian variety defined over a number field $k$.
Assume that the Mumford-Tate group of $A$ is the maximal possible, i.e.~the symplectic group 
$MT(A)=GSp_{2g}$ where $\dim A = g$,
and that the image of Galois action on the torsion has finite index in the group $GSp_{2g}(\Zl)$ 
of $\Zl$-points
of the symplectic group.

Then there are at most finitely many  $\Z$-form of the $l$-adic
cohomology theory $\Het(-\tensor \bar k,\Zl)$ restricted to  the category $\genn A k$, up to $\Aut(\kbar/k)$.
\end{example}

%
%
{\bf Proof} (sketch).
The Weil pairing corresponds 
 to the divisor corresponding to an ample line bundle over $A$,  
    and by compatibility with the cycle class map of a $\Z$-form and $\Hl=\Het(-\tensor \bar k,\Zl)$
    the non-degenerate Weil pairing
$$\omega : (\Hetl^1(A))^2 \lra \Hetl^0(A) = \Z_l$$
restricts to a pairing
 $$\omega : (H^1(A))^2 \lra H^0(A) \approx \Z,$$
which is easily seen to be non-degenerate.

    Now let $H_i$ be a $\Q$-form for $i=1,2$.

    Let $(x_i^1,...x_i^g,y_i^1,...,y_i^g)$ be a symplectic basis for $H_i^1(A)$.
    Then each is also a symplectic basis for $\Het^1(A)$,
    and so some $\sigma \in GSp(\Het^1(A),\omega)$ maps $H_1^1(A)$ to $H_2^1(A)$.
    The assumption on the Mumford-Tate group precisely means that such a
    $\sigma$ extends to $\sigma \in \Aut(\Hetl |_{\left<A\right>})$,
    and it follows from the fact that the cohomology of an Abelian variety is
    generated by $H^1$ that $\sigma(H_1) = H_2$.

 Finally, use the assumption on the Galois representation to see 
that there are at most finitely many $\Z$-forms.
\qed 

The proof above probably generalises to the following.

\begin{enonce}{Conjecture}[a generalisation of the example]
%

 Let $A$ be a motive of a smooth projective variety defined over a number field $k$. 
 Assume the Mumford-Tate group $G=Aut^\tensor(\Het(-\tensor\bar k, \Zl)_{|\left<A\right>})$
has the following property: 
 \begin{itemize}
 \item[] if $V_1$ and $V_2$ are abelian subgroups of 
$\Het(A\tensor\kbar,\Z_l)$
         which are both dense and of the same rank and 
         such that 
	 \begin{itemize}\item[]$GL(V_i) \cap G(\Z_l)$ is dense in $G(\Z_l)$ for i=1,2,
	 \end{itemize}
        then there is a $g \in G(\Z_l)$ such that $gV_1=V_2$ (setwise).
 \end{itemize}
  Then the conjectures  [2.1?,3.1?,3.2?,9.1?] of [Serre] imply that
 there are at most finitely many  $\Z$-forms of 
$\Het(-\tensor \bar k,\Zl)_{|\left<A\right>}$.
 \end{enonce}

Conjectures  [2.1?,3.1?,3.2?,9.1?] have analogues the cohomology theories
with coefficients in the ring of  finite adeles ${\Bbb A}^f$,
  cf. [Serre, 11.4?(ii), 11.5?], cf. also [Serre, 10.2?, 10.6?].

\subsection{\label{serre}Standard Conjectures and motivic Galois group}

Now we try to give a self-contained exposition of several conjectures 
on motivic Galois group which aapear related to our conjectures. 
Our exposition follows [Serre,\S1,\S3]

Let $k$ be a field of characteristic 0 which embeds into the field $\C$ of complex numbers; pick an embedding $\sigma:k\lra\C$. 

Assume Standard Conjectures and Hodge conjecture [Grothendieck, Kleiman]. 
Let $\Mot_k$ denote the category of pure motives over $k$ 
defined with the help of numerical equivalence of algebraic cycles
(or the homological equivalence, which should be the same 
by Standard Conjectures). $\Mot$ is a semi-simple category.

Let $E\in Ob \Mot$ be a motive; let $\left<E\right>$ denote
the least Tannakian subcategory of $\Mot_k$ containing $E$.

A choice of embedding $\sigma:k\lra\C$ defines 
an exact {\em fibre functor} 
$\Mot_k\lra\QVect$ corresponding to {\em the Betti realisation} 
$${\Hsing}_\sigma: \Mot_k\lra \QVect,\ \ \ E\mapsto \Hsing(E_\sigma(\C),\Q).$$ 

The scheme of automorphisms 
$MGal_{k,\sigma}=\Aut^\tensor({\Hsing}_\sigma:\Mot\lra\QVect)$ 
of the functor preserving the tensor product
is called {\em motivic Galois group of $k$}.
It is a linear proalgebraic group defined over $\Q$. 
Its category of $\Q$-linear representations is equivalent to $\Mot$. 
The group depends on the choice of $\sigma$. 

The motitivic Galois group of a motive $E$ is 
$\Aut^\tensor({\Hsing}_\sigma:\left<E\right>\lra\QVect)$. 

We now list several conjectures from [Serre]. 

\begin{enonce*}{Conjecture}[2.1?]
The group $\Aut^\tensor({\Hsing}_\sigma:\Mot_k\lra\QVect)$ is proreductive, i.e. a limit of liner reductive $\Q$-groups.
\end{enonce*}

Let $\boldsymbol 1$ denote the trivial morphism of rank 1, i.e. the cohomology 
of the point $Spec\, k$. 

\begin{enonce*}{Conjecture}[3.1?] The group
$\Aut^\tensor({\Hsing}_\sigma:\left<E\right>\lra\QVect)$ is the subgroup of
$GL({\Hsing}_\sigma(E))$ 
preserving the tensors corresponding to morphisms
 ${\boldsymbol 1}\lra E^{\tensor r} \tensor E^{\vee\tensor s}$, $r,s\geq 0$. 
\end{enonce*}

It is also conjectured that this group is reductive. 
Via the comparison isomorphism of \'etale and singular cohomology, 
$$\Hsing( E(\C), \Q ) \tensor \Q_l = \Het ( E\tensor\Qbar, \Q_l ),$$
the $\Ql$-points of $\Aut^\tensor({\Hsing}_\sigma:\left<E\right>\lra\QVect)(\Ql)$ act on the \'etale cohomology $\Het ( E\tensor\Qbar, \Q_l )$.
On the other hand, the Galois group $Gal(\Qbar/k)$ acts
on $\Qbar$ and therefore on $E\tensor\Qbar$. 
By functoriality, the Galois group acts by automorphisms of the functor of \'etale cohomology. 
 Hence, this gives rise to {\em $l$-adic representation associated to $E$}
  $$\rho_{k,l}: Gal(\Qbar/k)\lra \Aut^\tensor({\Hsing}_\sigma:\left<E\right>\lra\QVect)(\Ql).$$

\begin{enonce*}{Conjecture}[3.2?] Let $k$ be a number field. 
The image
of the $l$-adic representation associated with $E$
is dense in the group $\Aut^\tensor({\Hsing}_\sigma:\left<E\right>\lra\QVect)(\Ql)$ 
in the Zariski topology. 
\end{enonce*}

\begin{enonce*}{Conjecture}[9.1?] Let $k$ be a number field. 
The image $$
Im(\rho_{k,l}: 
Gal(\Qbar/k)\lra \Aut^\tensor({\Hsing}_\sigma:\left<E\right>\lra\QVect)(\Ql))$$ 
is open in $\Aut^\tensor({\Hsing}_\sigma:\left<E\right>\lra\QVect)(\Ql)$.
\end{enonce*}

\begin{enonce*}{Conjecture}[9.3?] Let $k$ be a number field. 
$ \Het ( E\tensor\Qbar, \Q_l )$ is semi-simple as a $Gal(\Qbar/k)$-module. 
\end{enonce*}

We suggest that the  conjectures [2.1?,3.1?,3.2?,9.1?,10.2?,10.3?.10.4?,10.7?,10.8?]
may be interpreted as saying there are only finitely
many $\Z$-forms of the \'etale cohomology 
$\Het(-,\Zl): \left<E\right>\tensor_k\Qbar\lra \QVect$,
up to Galois action. 
 There are similar conjectures for finite adeles instead of $\Q_l$,
  cf.~[Serre, 11.4?(ii), 11.5?], also [Serre, 10.2?, 10.6?].

\section{Speculations and remarks}

Standard conjectures claim there are algebraic cycles corresponding to various cohomological constructions. Model-theoretically it should mean
that something is definable in ACF and it is natural to expect
that such properties be useful in a proof of categoricity, i.e.~in the characterisation of the $\Q$-forms of \'etale cohomology theory. 

We wish to specifically point out the conjectures and properties 
involving smooth hyperplane sections, namely
{\em weak and strong Lefschetz theorems} and {\em Lefschetz Standard Conjecture}, cf.~[Kleiman,p.11,p.14]. 
Weak  Lefschetz theorem describes 
part of the cohomology ring of a smooth hyperplane section of a variety. Perhaps such a description can be useful in showing 
that a $\Q$-form extends uniquely to $Mot/K$ from the subcategory $\Mot/\Qbar$. 
An analogue of the weak Lefschetz theorem for the fundamental group was used in a similar way in [GavrDPhil, Lemma V.III.3.2.1],
see \ref{lefshetz:fund} for some details.
Namely, as is well-known, the fundamental group of a smooth hyperplane section of a smooth projective variety is  
essentially determined by the fundamental group of the variety. 
 [GavrDPhil, III.2.2] extends this to a somewhat technical weaker statement about arbitrary generic hyperplane sections.
An arbitrary variety can be represented as a generic hyperplane section of a variety defined over $\Qbar$ 
and this implies that, in some sense, the fundamental groupoid functor on the subcategory of varieties defined over $\Qbar$
``defines'' its extension to varieties defined over larger fields. 
The word ``defines'' is used in a meaning similar to model theoretic meaning of one first-order language definable in another.

\begin{question} Find a characterisation 
of the following families of functors:
 $$\Hsing( X(K_\tau), \Q ) : \Var/K \lra \Qvect,$$
  $$\Hsing( X(K_\tau), \C ) : \Var/K \lra \QHodge$$
where $\tau$ varies though isomorphisms of $K$ to $\C$,
or, almost equivalently, though locally compact locally 
connected topologies on $K$.
\end{question}

Note that Zilber [Zilber] {\em unconditionally} constructs a pseudo-exponential map $ex:\C^+\lra\C^*$ which satisfies the Schanuel conjecture. Of course, this map is not continuous (not even  measurable). 
 Hence we ask:
\begin{question} Construct a pseudo-singular cohomology theory 
which satisfies an analogue of the Schanuel conjecture and some other conjectures.
\end{question}

\subsection{Model theoretic conjectures}

Define model theoretic structures corresponding to the cohomology 
theories.

\begin{enonce}{Conjecture}
   The field is purely embedded into the structures
   corresponding to functors
\begin{itemize}
   \item[(i)]     $\Hsing:Var/\Qbar \lra \QVect$
   \item[(ii)]    $\Hsing(-,\Q):Var/\C \lra \QVect$
   \item[(iii)]   $\Hsing(-,\C):Var/\C \lra \QHodge$
\end{itemize}
Moreover, the structure (ii) is an elementary extension of (i) and
the cohomology ring $\Hsing(V,Q)$ is definable for every variety over C.
\end{enonce}

Several of the Standard Conjectures [Kleiman, \S4,p.11/9] claim that 
certain cohomological cycles (construction) correspond to algebraic 
cycles. This feels related to many of the conjectures above, 
in particular to the purity conjectures.

\begin{problem}
\begin{itemize}
\item[1.] Define a model-theoretic structure and language corresponding
 to the notion of a Weil cohomology theory, and formulate 
 a categoricity conjecture hopefully related to the Standard Conjectures ([Grothendieck, Kleiman])
 and conjectures on the motivic Galois Group and related Galois representations [Serre].

 \item[2.] Do the same in the language of functors, namely:
  \begin{itemize}
  \item[2.1.] Consider the family of cohomology theories on 
  $\Var/$K coming  from a choice of isomorphism $K\approx\C$.
  \item[2.2.] Define a notion of isomorphism of these/such cohomology theories,
    and  what it means to a "purely algebraic" property of such a theory.
   \item[3.3.] Find a characterisation of that family up to that notion of isomorphism
   by such properties. Or rather, show existance of such a characterisation is equivalent to a number of well-known conjectures
   such as the Standard Conjectures etc.
\end{itemize}
\end{itemize}
\end{problem}

\section{\label{pitop:uniq}Uniqueness properties of the topological fundamental groupoid functor of a complex algebraic variety}

\subsection{Statement of the conjectures}

Let $\VV$ be a category of varieties over a field $K$, let $\pi$
be a functor to groupoids such that $\Points \pi(X)=X(K)$
is the functor of $K$-points.
For $\sigma\in\Aut(K)$, define $\sigma(\pi)$ by
$$\Points\sigma(\pi)=\Points \pi(X) = X(K), \ \ \ 
\sigma(\Paths(x,y))=\Paths(\sigma(x),\sigma(y)),$$ 
$$ \ \ \ {\rm source}(\gamma)=\sigma({\rm source}(\gamma)), {\rm target}(\gamma)=\sigma({\rm target}(\gamma)),$$

For $K=\C$, an example of such a functor is the topological fundamental groupoid functor
$\pi_1^{top}(X(\C))$ of the topological space of complex points of an algebraic variety, 
and $\{\sigma(\pi_1^{top})\,:\,\sigma\in\Aut(\C)\}$ is the family of
all the topological fundamental groupoid functors associated with different choices
of a locally compact locally connected topology on $\C$. (Such a topology determines
a field automorphism, uniquely up to conjugation).

$Aut(K)$ acts by automorphisms of the source category, hence  all these 
(possibly non-equivalent!) 
 functors have the same properties in the language of functors, in particular
\begin{itemize}
\item[(0)] $\Ob\,\pi(X)=X(K)$ is the functor of $K$-points of an algebraic variety $X$
\item[(1)] preserve finite limits, i.e.~$\pi(X\times Y)=\pi(X)\times \pi(Y)$
\item[(2)] $\pi(X)$ is connected if $X$ is geometrically connected (i.e. the set of points $X(K)$ equipped 
with Zariski topology is a connected topological space)
\item[(3)] for $\tilde X\xra f X$ \'etale, the map $\pi(\tilde X)\xra{\pi(f)}\pi(X)$ of groupoids
		has the path lifting property of topological covering maps, namely
          \begin{itemize} 
 	  \item[] for $x =f(\tilde x),\tilde x\in \tilde X(K)$, for every path $\gamma\in \pi(X)$ starting at $x$,
                there is a unique path $\tilde\gamma\in \pi(\tilde X)$ such that 
                 ${\rm source}(\gamma)=\tilde x$ and $(\pi(f))(\tilde\gamma)=\gamma$.
          \end{itemize}
\end{itemize}



{\em A $\pi_1$-like functor} is a functor 
from a category of varieties to the category of groupoids
satisfying (0-3) above. 
Note that by (0) a $\pi_1$-functor comes 
equipped with a forgetful natural transformation to the functor of $K$-points.

\begin{enonce}{Conjecture}[$Z(\pi_1^{top})$] Each $\pi_1$-like functor on the category of smooth quasi-projective complex varieties
factors through the topological fundamental groupoid functor, up to a field automorphism. 

In detail: Let $\Var_\C$ be the category of smooth quasi-projective varieties over the field of complex numbers $\C$.
For each $\pi_1$-like functor $\pi:\Var_\C\lra\Groupoids$ there is a field automorphism $\sigma:\C\lra\C$ and a natural
transformation $\varepsilon:\pi_1^{top}\Lra\pi^\sigma$ such that the induced natural transformation $\Ob\pi_1^{top}\Lra\Ob\pi^\sigma$
on the functor of $\C$-points 
is identity.  
\end{enonce}

\begin{enonce}{Conjecture}[$Z(\pi_1,K)$] 
Let $K$ be an algebraically closed field.  Let $\Var_K$ be the category of smooth quasi-projective varieties over $K$.

There is a functor $\pi_1:\Var_K\lra\Groupoids$ such that 
for each $\pi_1$-like functor $\pi:\Var_K\lra\Groupoids$ there is a field automorphism $\sigma:K\lra K$ and a natural
transformation $\varepsilon:\pi_1\Lra \pi^\sigma$ such that the induced natural transformation $\Ob\pi_1\Lra\Ob\pi^\sigma$
on the functor of $K$-points 
is identity.  
\end{enonce}

\begin{remark} As stated, these conjectures are likely too optimistic. To get more plausible and manageable conjectures,
replace $\Var_K$ by a smaller category and add additional conditions on the $\pi_1$-like functors. 
The conclusion can also be weakened to claim there is a finite 
family of functors, rather than a single functor, through which $\pi_1$-like functors factor
up to field automorphism.

It may also be necessary to put extra structure on the fundamental groupoids. 
\end{remark}

\begin{remark} In model theory, it is more convenient to work with universal covering spaces 
rather than fundamental groupoids. Accordingly, model theoretic results are stated in the language
of universal covering spaces, sometimes with extra structure. 

 The conjectures above are motivated by questions and theorems about categoricity of certain structures. 
\end{remark}

\def\sepp{\mathrm{sep}}\def\Spec{\mathrm{Spec}\,}
\begin{remark} It is tempting to think that the right generalisation of the conjectures above should 
make use of the short exact sequence of \'etale fundamental groups (see [SGA1, XIII.4.3;XII.4.4])
$$1\lra 
\pi_1^{alg}(X\times_{\mathrm{Spec} k} \mathrm{Spec}k^\sepp, x)  \lra  \pi_1^{alg}(X, \bar x) 
\lra \pi_1(\Spec k, \mathrm{Spec}k^\sepp) =\Gal(k^\sepp/k)\lra 1$$
where $X$ is a scheme over a field $k$, $k^\sepp$ is a separable closure of $k$, and  $x:\Spec k^\sepp \lra X\times_{\mathrm{Spec} k} \mathrm{Spec}k^\sepp$ is a geometric point of $X\times_{\mathrm{Spec} k} \mathrm{Spec}k^\sepp$,
and $\bar x : \Spec k^\sepp \lra X\times_{\mathrm{Spec} k} \mathrm{Spec}k^\sepp\lra X $ is the corresponding geometric point of $X$.
 
In fact such a sequence could be associated with a morphism $X\lra S$ 
admitting a section and satisfying certain assumptions [SGA 1, XIII.4].

These short exact sequences comes from pullback squares
\def\pullbacksquare#1#2#3#4#5#6{\xymatrix{ {#1} \ar[r]|{} \ar@{->}[d]|{#2} & {#4} \ar[d]|{#5} \\ {#3}  \ar[r] & {#6} }}
$$\pullbacksquare{X\times_{\mathrm{Spec} k} \mathrm{Spec}k^\sepp}{}{X}{\Spec k^\sepp}{}{\Spec k}
\ \ \ \ \ \ \ \ \ \ \ \ \  \ \ \ 
\pullbacksquare{X_s}{}{X}{s}{}{S}
$$

\end{remark}

We find the following conjecture plausible and hope its statement clarifies the arithmetic nature of our conjectures.
It is perhaps the simplest conjecture not amendable to model theoretic analysis because it uses bundles. 
In the next subsection we list several partial positive results.

For a variety $X$, let $\genn X K$ denote the category  whose objects are the finite Cartesian powers 
of $X$, and morphisms are morphisms of algebraic varieties defined over $K$; 
we let $X^0$ to be a variety consisting of a single $K$-rational point.

\begin{enonce}{Conjecture}[$Z(\pi_1,L_A^*)$] Let $K$ be an algebraically closed field of zero characteristic, 
$A$ an Abelian variety defined over a number field $k$. Let $L_A$ be an ample line bundle over $A$
 and $L_A^*$ be the corresponding ${\Bbb G}_m$-bundle. 
Further assume that the Mumford-Tate group of $A$ is the maximal possible, i.e.~the general symplectic group,
$$MT(A)=\GSpZ$$
and that the image of Galois action on the torsion has finite index in the group of ${\Bbb {\hat{Z}}}$-points
of the symplectic group.

Then there is a finite family of $\pi_1$-like functors $\Pi_1$ such that each $\pi_1$-like functor 
 on the full subcategory $\genn {L_A^*} K$ 
consisting of the Cartesian powers of the ${\Bbb G}_m$-bundle $L_A^*$,
factors via an element of $\Pi_1$.
\end{enonce}

These functors in $\Pi_1$ correspond to different embeddings of the field of definition of $A$ into the field of complex numbers. 

The following conjecture is probably within reach, at least if we replace the fundamental groupoid functor by its residually finite part. 

Model theoretic methods of $[\mathrm{BH^2K^2 14}]$ probably allow to replace $\C$ by a countable algebraically closed subfield. 
Methods of [GavrDPhil,III.5.4.7], 
cf.~\S\ref{lefshetz:fund},  
probably reduce the remaining part of the conjecture 
to properties of complex analytic topology
and normalisation of varieties.

\begin{enonce}{Conjecture}[$Z(\pi_1,\Qbar\subset \C)$]
Let $\Var_\C$ be the category of smooth quasi-projective varieties over $\C$,
and let $\Var_{\Qbar}$ be its category consisting of varieties and morphisms defined over $\Qbar$.

Assume that $\pi:\Var_\C\lra \Groupoids$ is a $\pi_1$-like functor 
which coincides with the topological fundamental groupoid 
$\pi_1^{top}$ for varieties and morphisms defined over $\Qbar$, i.e. 
for each variety $V$ in $\Var_{\Qbar}$ 
and each morphism $V\xra f W$ in $\Var_{\Qbar}$ 
 it holds
$\pi(V)=\pi_1^{top}(V)$ and $\pi(f)=\pi_1^{top}(f)$. 

Then there exist a field automorphism $\sigma\in \Aut(\C/\Qbar)$
such that $\pi\circ\sigma$ and $\pi_1^{top}$ are equivalent.
\end{enonce}

There are a number of theorems and conjectures which can be seen as saying that, 
up to finite index, Galois action is described by geometric, algebraic or topological structures;
our conjectures can also be seen in this way.

%

\subsection{\label{results}Partial positive results} These conjectures are closely related to {\em categoricity} 
theorems in model theory, and this led to several partial positive results 
about the full subcategories $\gen {K^*}$ of algebraic tori in arbitrary characteristic,
$\gen E$ powers of an elliptic curve over a number field, a weaker result about
$\gen A$ powers of an Abelian variety over a number field, a still  weaker result
about $\gen V$ powers of a smooth projective variety whose fundamental group 
satisfies a group theoretic property of being subgroup separable,
a strengthening of residually finite. 

Note that the first three categories are linear in the sense that the the groups $\Aut(K^*)$,
$\Aut_{End\,E-mod}(E(K))$, and $\Aut_{End\,A-mod}A(K)$ act on the set of $\pi_1$-like functors on 
the respective categories $\gen {K^*}$, $\genn E K$, and $\genn A K$. This is so because
these groups act on these categories.

Below, we list several known results, translated from categoricity theorems available in model theory 
literature. We list the corresponding category $\VV$ and additional properties requires 
of the functors in the family.

\begin{enumerate}
\item \label{ref:Gm} [BaysZilber,Th.2.1] $\char K=0$,  $\gen {K^*}$ 

\item\ [BaysDPhil, Th.4.4.1; GavrK, Prop.2]  $\char K=0$, $\genn E K$ where $E$ is an elliptic curve defined over a number field $k$
with a $k$-rational point $0\in E(k)$; 
$\pi(E,0,0)\approx \Z^2$; a finite family

\item\  [BaysZilber,Th.2.2] $\char K=p>0$, $\gen {K^*}$; 
there is a $\pi_1$-like functor $\pi_1:\gen {K^*}\lra \Groupoids$ such that each $\pi_1$-like functor
$\pi:\gen {K^*}\lra \Groupoids$ factors via $\pi_1:\gen {K^*}\lra \Groupoids$ up to $\Aut(K/\Fpbar)$
provided
\begin{itemize}
 \item $\pi(K^*,1,1)\approx {\Bbb Z}[1/p]$  
\item 
the restrictions  ${\pi_1}_{|\bar{\Bbb F}_p}$ and $\pi_{|\bar{\Bbb F}_p}$ 
to $\Fpbar$-rational points coincide:
$$\pi_{|\bar{\Bbb F}_p}=\pi'_{|\bar{\Bbb F}_p}$$
 \end{itemize}

\item\ [BaysDPhil, Th.4.4.1] \label{ab:cat} $\char K=0$, $\genn A K$
 where $A$ is an Abelian variety defined over a number field $k$ with a $k$-rational point $0\in A(k)$, 
 \begin{itemize}
 \item $\pi(A,0,0)={\Bbb Z}^{2\dim A}$
 \item  for any two functors $\pi,\pi'$ in $\mathcal F$,  
the corresponding fundamental group functors coincide
$$\pi(A,0,0)=\pi'(A,0,0)$$
and further, for $p:\tilde A\lra A$ is \'etale, $\gamma\in\pi(A,0,0)=\pi'(A,0,0)$,
$\pi(\gamma_\pi)=\pi'(\gamma_{\pi'})=\gamma$, it holds that
$${\rm source}(\gamma_\pi)={\rm source}(\gamma_{\pi'})\ \ \ {\rm 
implies }\ \ \ \ {\rm target}(\gamma_\pi)={\rm target}(\gamma_{\pi'})$$ 
 \end{itemize}

\item\ [GavrDPhil,III.5.4.7] \label{w-cat} $\char K=0$, ${\rm card\,}K=\aleph_1$, and $\genn V K$
 where $V$ is an smooth projective variety defined over a number field $k$
with a $k$-rational point $0\in V(k)$ such that the universal covering space of $V(\C)$ is holomorphically complex, for some embedding $K\hookrightarrow \C$, 
and  its fundamental groups $\pi_1(V(\C),0,0)^n$ are subgroup separable for each $n>0$;
recall a group $G$ is subgroup separable iff for each finitely generated subgroup $H<G$
and $h\not\in H$ there is a morphism $f:G\lra G_H$ into a finite group $G_H$ such that
$h\not\in f(H)$. 

there is a $\pi_1$-like functor $\pi_1:\gen V\lra \Groupoids$ such that a $\pi_1$-like functor
$\pi:\gen {V}\lra \Groupoids$ factors via $\pi_1:\gen {V}\lra \Groupoids$ up to $\Aut(K/\Qbar)$
provided
\begin{itemize}
 \item $\pi(V,0,0)\approx \pi_1^{top}(V(\C),0,0)$ 
\item 
the restrictions  ${\pi_1}_{|\bar{\Bbb Q}}$ and $\pi_{|\bar{\Bbb Q}}$ 
to $\Qbar$-rational points coincide:
$$\pi_{|\bar{\Bbb Q}}=\pi'_{|\bar{\Bbb Q}}$$
 \end{itemize}


\item We wish to mention the work of [HarrisDPhil, DawHarris] on Shimura curves, which does not quite fit 
in our framework. To interpret their results, one needs to consider $\pi_1^{top}$ 
as a functor to groupoids {\em with extra structure}.  

\end{enumerate}

Conjectures on independence of Galois representations of non-isogenious curves probably imply our conjectures 
for the full subcategory
$\genn {E_1\times\ldots\times E_n} K$
generated by a finite product of elliptic curves $E_1,\ldots,E_n$ over a number field.

Consider the family of $\pi_1$-like functors with Abelian fundamental groups.
this requires weakening of the uniqueness in the path-lifting property (3). 
Is it easier to prove that each such functor factors via $\pi_1^{top}$ up to a  field automorphisms?


\subsection{Mathematical meaning of the conjectures. Elements of proof of the conjectures}

Here we try to explain the arithmetic and geometric meaning of the conjectures. 
In a sense, the conjectures say that $\Gal(\Qbar/\Q)$ and $\Aut(K/\Q)$ are large enough.
We try to show below in what sense, by showing possible obstructions/difficulties in proof.

\subsubsection{Galois action on roots of unity and Kummer theory}
Consider the infinite sequence  $exp(2\pi i/n)$ of roots of unity. This sequence can be obtained
topologically: take the loop $\gamma$ generating $\pi(\C^*,1,1)\approx \Z$, the \'etale 
morphism $z^n:\C^*\lra\C^*$ and lift $\gamma$ uniquely to a path $\tilde\gamma_n$ starting at $1\in \C^*$. 
Then $exp(2\pi/n)$ is the end-point of $\gamma_n$. This construction shows that a $\pi_1$-like functor on 
the category $\gen {K^*}$ determines a distinguished sequence $\xi_n,n\geq 0$ of roots of unity. 
 Hence, our conjectures require that the Galois group acts transitively on the set of sequences of roots of unity
associated with $\pi_1$-like functors. 

Consider a $\pi_1$-like functor on the category $\gen {K^*}$. As noted above, group automorphisms $\Aut(K^*)$ 
of $K^*$ act on the set of these functors. Hence, item (\ref{ref:Gm}) requires that multiplicative group automorphisms 
$\Aut(K^*)$ and field automorphisms $\Aut(K/\Q)$ have the same orbits on the sequences $\xi_n,(\xi_{mn})^{m}=\xi_n,m,n>0$ 
of roots of unity.  

Kummer theory arises in a similar way if we consider endpoints of liftings of paths joining $1$ and 
arbitrary elements $a_1,...,a_n$.

\subsubsection{Elliptic curves and Abelian varieties. Kummer theory and Serre's open image theorem
for elliptic curves.}
Kummer theory for elliptic curves and Abelian varieties arises in the same way if we consider 
$\pi_1$-like functors on the category $\gen A$ generated by an Abelian variety. 

 Similarly, our conjectures 
about $\pi_1$-like functors on $\genn A K$ require that the action of $Aut_{EndA{\rm-mod}}(A(K))$
and $\Gal(\Qbar/k)$ on the torsion points do not differ much. This is true for elliptic curves
but fails for Abelian varieties of $\dim A>1$, hence the extra assumption in (4) on the family of $\pi$-like functors.

\subsubsection{\label{lefshetz:fund}Arbitrary variety. Etale topology and an analogue of Lefshetz theorem for the fundamental group}

To prove item (5), we need several facts about \'etale topology. Most of these facts are well-known for smooth varieties;
what we use is that they hold ``up to finite index'' for arbitrary (not necessarily smooth or normal) subvarieties
of a smooth projective variety. 

\def\limVC{\varprojlim\tilde V(\C)}
Consider the inverse limit $\varprojlim \tilde V(\C)$ of finite \'etale covers $\tilde V(\C) \lra V(\C)$  of a complex algebraic 
variety $V$. The universal analytic covering map $U\lra V(\C)$ gives rise to covering maps $U\lra \tilde V(\C)$ and hence
a map $U\lra \limVC$. Zariski topology on the \'etale covers makes $\limVC$ into a topological space. Hence 
there are two topologies on $U$ -- the complex analytic topology and the ``more algebraic"  topology
on $U$ induced from the map $U\lra\varprojlim V(\C)$. Call the latter {\em \'etale} topology on $U$.

To prove item (\ref{w-cat}) we use that these two topologies are similar and  nicely related. In particular,
\begin{itemize}
\item Closed irreducible sets in \'etale topology are closed irreducible in complex analytic topology (by definition).
\item For a set closed in \'etale topology, its irreducible components in complex analytic topology 
  are also closed in \'etale topology [GavrDPhil,III.1.4.1(4,5)].
\item The image of an \'etale closed irreducible subset of $U\times\ldots\times U$ under a coordinate projection is \'etale closed [GavrDPhil, III.2.2.1].
\end{itemize} 

Note that this is easy to see that {\em connected} components of a set closed in \'etale topology are also 
closed in \'etale topology, 
and hence that the properties above holds for smooth or normal closed sets. 

Let $f:W\lra V$ be a morphism of varieties, and let $f_*:U_W\lra U_V$ be the map of the universal covering 
spaces of $W(\C)$ and $V(\C)$.  We may assume that $V$ is smooth projective but it is essential that $W$ is
arbitrary. In applications, $W$ is an arbitrary closed subvariety of a Cartesian power of a fixed variety $V$. 
\begin{itemize}
\item 
       if $f:W\lra V$ is proper, then the image $f(U_W)$ is closed in $U_V$ in \'etale topology
\end{itemize} 
This is related to the following geometric fact [GavrDPhil, V.3.3.6, V.3.4.1]:
\begin{itemize}
\item If $f:W(\C)\lra V(\C)$ is a morphism of smooth normal algebraic varieties, $g$ a generic point of $V(\C)$
and $W_g=f^{-1}(g)$ then 
$$
\pi_1(W_g,w,w)\lra\pi_1(W,w,w)\lra \pi_1(V,g,g)$$
is exact up to finite index
\item Moreover, if $f(W(\C))$ is dense in $V(\C)$, then $\pi_1(W,w,w)\lra \pi_1(V,g,g)$ is surjective.
\end{itemize} 
In fact we use a generalisation of this, namely that it holds up to finite index for arbitrary varieties
if one considers the image of the fundamental group in the ambient smooth projective variety.


{\bf Acknowledgement.} Ideas and proofs were strongly influenced by extensive conversations with Martin Bays.
I thank A.Luzgarev and V.Sosnilo for useful discussions. 
I thank Sergei Sinchuk for helpful discussions.
I also thank Maxim Leyenson for comments on a late draft. 
I thank Yves Andre for several corrections. 

 Support from Basic Research Program of the National Research University Higher                                                                             
School of Economics is gratefully acknowledged. This study was partially supported by                                                                          
the grant 16-01-00124-a of Russian Foundation for Basic Research.

\end{document}